\documentclass{amsart}
\usepackage{epsfig}
\usepackage{amsmath}
\usepackage{amsfonts}
\begin{document}

\newtheorem{thm}{Theorem}
\newtheorem{lem}[thm]{Lemma}
\newtheorem{cor}[thm]{Corollary}
\newtheorem{conj}[thm]{Conjecture}
\newtheorem{qn}{Question}
\newtheorem{pro}{Proposition}[section]

 \theoremstyle{definition}
\newtheorem{defn}{Definition}[section]

\theoremstyle{remark}
\newtheorem{rmk}{Remark}[section]

\def\square{\hfill${\vcenter{\vbox{\hrule height.4pt \hbox{\vrule width.4pt
height7pt \kern7pt \vrule width.4pt} \hrule height.4pt}}}$}

\def\Z{\Bbb Z}
\def\R{\Bbb R}
\newenvironment{pf}{{\it Proof:}\quad}{\square \vskip 12pt}

\title[complex Kleinian groups]{The complex Kleinian groups with an invariant totally geodesic submanifold}
\author{ Baohua Xie}

\address{College of Mathematics and Econometrics \\ Hunan University \\ Changsha, 410082, China}
\email{xiexbh@gmail.com}
 \maketitle

\begin{abstract}
In this paper, we characterize discrete subgroups of $\mathbf{PU}(2,1)$, holomorphic isometric group of complex hyperbolic space,
which keep invariant an invariant totally geodesic submanifold.
 \end{abstract}

\section{Introduction}
A Kleinian group $G$ is Fuchsian if it keeps invariant some circular disc $\mathbf{U}$. We can regard $\mathbf{U}$
as being hyperbolic plane $\mathbf{H}^{2}$, so that a Fuchsian group is a discrete subgroup of $\mathbf{PSL}(2, \mathbf{R})$.
The following theorem can be found in V.G.18 in Maskit's book \cite{m}.

 \begin{thm}Let $G\subset \mathbf{SL}(2, \mathbf{C})$ be a non-elementary Kleinian group in which
 $tr^{2}(g)\geq 0$ for all $g\in G$. Then $G$ is Fuchisn.
\end{thm}
This theorem gave a characterization of those Kleinian groups that are Fuchsian.

In this paper we are interested  in discrete subgroups acting on complex hyperbolic space
which have an invariant totally geodesic submanifold with codimension $2$. There are no geodesic hypersurface
in complex hyperbolic space. But there are two kinds of totally geodesic submanifolds: Complex line and
Langrangian  plane. In section 3, we will characterize complex Kleinian groups which keep invariant  an invariant totally geodesic submanifold .
Recently, J. Kim \cite{kim} proved a similar result
in quaternionic hyperbolic case.

\section{Complex hyperbolic space}
\subsection{Siegel domain}
Let $\mathbf{C}^{2,1}$ be a complex vector space of dimension 3 with a Hermitian form of Sigature
$(2,1)$ given by $$\langle \mathbf{z}, \mathbf{w}\rangle=z_{1}\overline{w_{3}}+z_{2}\overline{w_{2}}+z_{3}\overline{w_{1}}.$$
We consider the subspaces
$$V_{-}=\{\mathbf{z}\in \mathbf{C}^{2,1}: \langle \mathbf{z}, \mathbf{z}\rangle<0 \},$$
$$V_{0}=\{\mathbf{z}\in \mathbf{C}^{2,1}: \langle \mathbf{z}, \mathbf{z}\rangle=0 \}.$$

Let $\mathbf{P}: \mathbf{C}^{3}-\{0\}\longrightarrow \mathbf{C}\mathbf{P}^{2}$ be the canonical
projection. Then complex hyperbolic space is defined to be $\mathbf{H}_{\mathbf{C}}^{2}=\mathbf{P}(V_{-})$ and $\partial\mathbf{H}_{\mathbf{C}}^{2}=\mathbf{P}(V_{0})$ is its boundary. Using non-homogeneous coordinates we can write
$\mathbf{H}_{\mathbf{C}}^{2}$ as

$$\mathbf{H}_{\mathbf{C}}^{2}=\{(z_{1},z_{2})\in \mathbf{C}^{2}: 2\Re z_{1}+|z_{2}|^{2}<0\}$$
and also, for $\partial\mathbf{H}_{\mathbf{C}}^{2}$ we have
$$\partial\mathbf{H}_{\mathbf{C}}^{2}=\{(z_{1},z_{2})\in \mathbf{C}^{2}: 2\Re z_{1}+|z_{2}|^{2}=0\}.$$

Given a point $z$ of $\mathbf{C}^{2}\subset \mathbf{C}\mathbf{P}^{2}$ we may lift $z=(z_{1},z_{2})$ to a point $\mathbf{z}$ in $\mathbf{C}^{2,1}$,
called the standard lift of $z$, by writing $\mathbf{z}$ in non-homogeneous coordinates as

$$\mathbf{z}=\begin{bmatrix}z_{1}\\
z_{2}\\
1\end{bmatrix}.$$

Specially, the standard lifts of $\mathbf{0}=(0,0)$  and $\infty$ are as follows
$$\mathbf{0}=\begin{bmatrix}0\\
0\\
1\end{bmatrix}, \ \ \infty=\begin{bmatrix}1\\
0\\
0\end{bmatrix}.$$
Complex hyperbolic space is a $2$-complex dimensional complex manifolds. The Bergman metric on $\mathbf{H}_{\mathbf{C}}^{2}$
is defined by the distance function given by the formula
\[\cosh^{2}\big(\frac{\rho(z,w)}{2}\big) =\frac{\langle \mathbf{z},\mathbf{w}\rangle \langle \mathbf{w},\mathbf{z}\rangle}{\langle \mathbf{z},\mathbf{z}\rangle\langle \mathbf{w},\mathbf{w}\rangle}. \]

\subsection{Complex hyperbolic isometries}

The group of biholomorphic transformations of $\mathbf{H}_{\mathbf{C}}^{2}$ is  $\mathbf{PU}(2,1)$, the projection of the Unitary
group $\mathbf{U}(2,1)$  preserving the Hermitian form given in 2.1. In this work we prefer to consider instead the  group $\mathbf{SU}(2,1)$
of matrices which are unitary with respect to $\langle \cdot,\cdot\rangle$,
and have determinant 1. The group $\mathbf{SU}(2,1)$ is a 3-fold covering of $\mathbf{PU}(2,1)$, a direct analogue of the fact that $\mathbf{SL}(2,\mathbf{C})$ is the double cover of $\mathbf{PSL}(2,\mathbf{C})$.

The general form of an element of $A\in \mathbf{SU}(2,1)$ and its inverse are

$$A= \begin{bmatrix} a& b& c\\
d& e& f\\
g &h&
j\end{bmatrix}, \ \ A^{-1}= \begin{bmatrix} \overline{j}& \overline{f}& \overline{c}\\
\overline{h}& \overline{e}& \overline{b}\\
\overline{g} &\overline{d}&
\overline{a}\end{bmatrix}. $$

As the composition of an element of $\mathbf{SU}(2,1)$ with its inverse is the  identity, we
obtain a list of equations that the matrix entries in $A$ must satisfy
\[a\overline{j}+b\overline{h}+c\overline{g}=1,a\overline{f}+b\overline{e}+c\overline{d}=0,
a\overline{c}+b\overline{b}+c\overline{a}=0,d\overline{j}+e\overline{h}+f\overline{g}=0,\]
\[d\overline{f}+e\overline{e}+f\overline{d}=1,
g\overline{j}+h\overline{h}+j\overline{g}=0,
a\overline{j}+d\overline{f}+g\overline{c}=1,b\overline{j}+e\overline{f}+h\overline{c}=0,\]
\[c\overline{j}+f\overline{f}+j\overline{c}=0,
a\overline{h}+d\overline{e}+g\overline{b}=0,b\overline{h}+e\overline{e}+h\overline{b}=1,
a\overline{g}+d\overline{d}+g\overline{a}=0.\]

There exist three kinds of holomorphic isometries of $\mathbf{H}_{\mathbf{C}}^{2}$.

{\rm(i)} Loxodromic isometries, each of which fixes exactly two points of $\partial\mathbf{H}_{\mathbf{C}}^{2}$.
One of these points is attracting and the other repelling.

{\rm(ii)} Parabolic isometries, each of which fixes exactly one point of $\partial\mathbf{H}_{\mathbf{C}}^{2}$.

{\rm(iii)} Elliptic isometries, each of which fixes at least one point of $\mathbf{H}_{\mathbf{C}}^{2}$.

\subsection{Geodesic submanifolds}
Unlike the real hyperbolic space case, there are no totally geodesic submanifolds of codimension 1 in complex hyperbolic space.
But there are two kinds of totally geodesic 2-dimensional subspaces of complex hyperbolic space.
Namely:

{\rm (i)} Complex lines $L$, which have constant curvature $-1$, and

{\rm(ii)} totally real Langrangian planes $R$, which have constant curvature $-\frac{1}{4}$.

Every Complex line $L$ is the image under some element of $\mathbf{SU}(2,1)$ of the
complex line $L_{1}$ with polar vector $n_{1}=(0,1,0)^{t}$. The complex line  $L_{1}$ has the following
form
$$L_{1}=\{(z_{1},z_{2})^{t}\in \mathbf{H}_{\mathbf{C}}^{2}: z_{2}=0\}.$$

The subgroup of $\mathbf{SU}(2,1)$ stabilizing  $L_{1}$ is thus conjugate to the group $S(\mathbf{U}(1)\times \mathbf{U}(1,1))<\mathbf{SU}(2,1)$. The stabilizer
of every other Complex line is conjugate to this subgroup. Every Lagrangian plane is the image under some element
of $\mathbf{SU}(2,1)$ of the standard real Lagrangian plane $L_{\mathbf{R}}$, where both coordinates are real:

$$L_{\mathbf{R}}=\{(z_{1},z_{2})^{t}\in \mathbf{H}_{\mathbf{C}}^{2}: \Im(z_{1})=\Im(z_{2})=0\}.$$

A complex Kleinian group $G\subset\mathbf{SU}(2,1)$ is $\mathbf{R}-$Fuchsian (or $\mathbf{C}-$Fuchsian)
if it keeps invariant some Langrangian plane ( or Complex line).

\subsection{Cartan's angular invariant}

Let $z_{1}, z_{2}, z_{3}$ be three distinct points of $\partial\mathbf{H}_{\mathbf{C}}^{2}$ with lifts
$\mathbf{z}_{1}, \mathbf{z}_{2},\mathbf{ z}_{3}$. Cartan's angular invariant is defined as follows:
$$\mathbb{A}(z_{1}, z_{2}, z_{3})=\arg(-\langle \mathbf{z}_{1}, \mathbf{z}_{2} \rangle)\langle \mathbf{z}_{2}, \mathbf{z}_{3} \rangle)\langle \mathbf{z}_{3}, \mathbf{z}_{1} \rangle).$$
The angular invariant is independent of the chosen lift $\mathbf{z}_{j}$. It is clear that applying an element of $\mathbf{SU}(2,1)$
to our triple of points does not change the Cartan invariant. The properties of $\mathbb{A}$ may be found  in Goldman's book
\cite{g}. In the next proposition we highlight some of them.
\begin{pro}
Let $z_{1}, z_{2}, z_{3}$ be three distinct points of $\partial\mathbf{H}_{\mathbf{C}}^{2}$ and $\mathbb{A}(z_{1}, z_{2}, z_{3})$ be
their angular invariant. Then

{\rm (i)} $\mathbb{A}\in [-\frac{\pi}{2},\frac{\pi}{2}]$;

{\rm (ii)} $\mathbb{A}=\pm \frac{\pi}{2}$ if and only if  $z_{1}, z_{2}, z_{3}$ all lie on a Complex line;

{\rm (ii)} $\mathbb{A}=0$ if and only if  $z_{1}, z_{2}, z_{3}$ all lie on Langrangian plane;
\end{pro}

\subsection{The Kor\'{a}nyi-Reimann cross-ratio}
Cross-ratios were introduced to complex hyperbolic space by Kor\'{a}nyi and Reimann\cite{km}. We suppose that $z_{1},z_{2},z_{3},z_{4}$ are four
distinct points of $\partial\mathbf{H}_{\mathbf{C}}^{2}$. Let $\mathbf{z}_{1},\mathbf{z}_{2},\mathbf{z}_{3},\mathbf{z}_{4}$ be corresponding lifts
in $V_{0}\subset \mathbf{C}^{2,1}$. The Kor\'{a}nyi-Reimann cross-ratio of this four points is defined to be
$$\mathbb{X}=\frac{\langle \mathbf{z_{3}},\mathbf{z_{1}}\rangle \langle \mathbf{z_{4}},\mathbf{z_{2}}\rangle}{\langle \mathbf{z_{4}},\mathbf{z_{1}}\rangle\langle \mathbf{z_{3}},\mathbf{z_{2}}\rangle}.$$
$\mathbb{X}$ is invariant under $\mathbf{SU}(2,1)$ and independent of the chosen lifts.

In order to study the configure space of quadruple of points $z_{1},z_{2},z_{3},z_{4}$ in the boundary   $\partial\mathbf{H}_{\mathbf{C}}^{2}$, J.R. Parker and I.D. Platis\cite{pp}defined
other cross-ratios by choosing different ordering of the four points. Given distinct points $z_{1},z_{2},z_{3},z_{4}\in \partial\mathbf{H}_{\mathbf{C}}^{2}$, they defined
$$\mathbb{X}_{1}=[z_{1},z_{2},z_{3},z_{4}],\mathbb{X}_{2}=[z_{1},z_{3},z_{2},z_{4}],\mathbb{X}_{3}=[z_{2},z_{3},z_{1},z_{4}].$$

Moreover, they showed that all three of $\mathbb{X}_{1},\mathbb{X}_{2}$ and $\mathbb{X}_{3}$ are real if and only if the four points either lie in the same Complex line or on the same Lagrangian plane. That is,

\begin{pro}\cite{pp} Suppose that $\mathbb{X}_{1},\mathbb{X}_{2}$ and $\mathbb{X}_{3}$ are all real.

(i) If $\mathbb{X}_{3}=-\mathbb{X}_{2}/\mathbb{X}_{1}$ then points $z_{j}$ all lie on a Complex line.

(ii) If $\mathbb{X}_{3}=\mathbb{X}_{2}/\mathbb{X}_{1}$ then points $z_{j}$ all lie on a Lagrangian plane.
\end{pro}

\section{Main results}

In this section, we prove our main theorem. Before stating our results,  we begin by recalling some notions.
 A subgroup $G\subset \mathbf{SU}(2,1)$  is called elementary if it has a finite $G$-orbit in $\overline{\mathbf{H}_{\mathbf{C}}^{2}}$. Otherwise we call $G$ a non-elementary group. We say $G$ is Fuchsian if $G$ is either $\mathbf{C}$-Fuchsian
or $\mathbf{R}$-Fuchsian. Then the statement of our  result is almost the same as real hyperbolic case by Maskit.

 \begin{thm} Let $G$ in $\mathbf{SU}(2,1)$ be a non-elementary complex hyperbolic Kleinian
group. If the trace of every element of $G$ is real then $G$ is Fuchsian.
\end{thm}
\begin{pf}
If $G$ is $\mathbf{R}$-Fuchsian, then $G$ is conjugation  to a subgroup of $\mathbf{SO}(2,1)$. If $G$ is $\mathbf{C}$-Fuchsian, then $G$ is conjugation to a subgroup of $S(\mathbf{U}(1)\times \mathbf{U}(1,1))<\mathbf{SU}(2,1)$.
So every element of $G$ has real trace in both case.

 We need to prove the converse. Since $G$ is a non-elementary group, there exist a loxodromic element $A\in G$ with real trace. Because $\mathbf{PU}(2,1)$ acts two points homogeneously on $\partial\mathbf{H}_{\mathbf{C}}^{2}$ and the trace is invariant under conjugation, we may assume that $A$ fixes $0$ and $\infty$.
Now select any $B$ in $G$. In term of matrices we can write
$$A= \begin{bmatrix} t& 0& 0\\
0& 1& 0\\0&0& \frac{1}{t}\end{bmatrix}
,\ B= \begin{bmatrix} a& b& c\\
d& e& f\\g&h& j\end{bmatrix},$$
where each matrix is in $\mathbf{SU}(2,1)$.  We can assume that $t>1$.

Next, write
\begin{eqnarray*}
\begin {array}{ll}
t_{1}=tr(AB)=ta+e+\frac{j}{t}\in \mathbf{R}
\\t_{2}=tr(A^{-1}B)= \frac{a}{t}+e+jt\in \mathbf{R}\\
t_{3}=tr(B)=a+e+j\in \mathbf{R}
\end{array}
\end{eqnarray*}

Since $A$, $AB$ and $A^{-1}B$ have real trace, $t_{1}$, $t_{2}$ and $t_{3}$ are real. Solving for $a, e, j$, we find that $a, e, j \in \mathbf{R}$. This shows that every element of $G$ has real diagonal elements.

Let $B'=BAB^{-1}$. Then the matrix of $B'$ is

$$B'= \begin{bmatrix} ta\overline{j}+b\overline{h}+c\overline{g}/t& ta\overline{f}+b\overline{e}+c\overline{d}/t& ta\overline{c}+b\overline{b}+c\overline{a}/t\\
td\overline{j}+e\overline{h}+f\overline{g}/t& td\overline{f}+e\overline{e}+f\overline{d}/t& td\overline{c}+e\overline{b}+f\overline{a}/t\\tg\overline{j}+h\overline{h}+j\overline{g}/t&tg\overline{f}+h\overline{e}+j\overline{d}/t& tg\overline{c}+h\overline{b}+j\overline{a}/t\end{bmatrix}.$$

The fixed points of $B'$ are $B(0)$ and $B(\infty)$. In homogeneous coordinate,

$$B(0)= \begin{bmatrix} c\\
f\\j\end{bmatrix}
, \ B(\infty)= \begin{bmatrix} a\\
d\\g\end{bmatrix}.$$

To begin with, we prove that $c$ is a real number or a pure imaginary number.

First, note that
$$B^{2}= \begin{bmatrix} a^{2}+bd+cg& *& *\\
*& db+e^{2}+fh& *\\*&*& cg+hf+j^{2}\end{bmatrix}$$ also has real diagonal elements.
So $a^{2}+bd+cg,db+e^{2}+fh, cg+hf+j^{2}$ are real. Then $cg$ is real.

Since the diagonal elements of $B'$ are real, we have $ta\overline{j}+b\overline{h}+c\overline{g}/t$,
$tg\overline{c}+h\overline{b}+j\overline{a}/t$ are real. Thus $ta\overline{j}+b\overline{h}+c\overline{g}/t+
tg\overline{c}+h\overline{b}+j\overline{a}/t$ is real. Hence, $tg\overline{c}+c\overline{g}/t$ is real.

Suppose  $c=r_{1}e^{i\theta},g=r_{2}e^{-i\theta}$,$r_{1},r_{2}\neq 0$.
Then we have
$$r_{1}r_{2}\sin(2\theta)/t-tr_{1}r_{2}\sin(2\theta)=0.$$ That is,

$$(r_{1}r_{2}/t-tr_{1}r_{2})\sin(2\theta)=0.$$ Then $\sin(\theta)=0$ or $\cos(\theta)=0$.
Therefore $c$ is real or a pure imaginary number.

Now we claim that $c\neq 0$.

If $c=0$, then $f=0$ by the identity $c\overline{j}+f\overline{f}+j\overline{c}=0$.
 So $B(0)=0$. But $A$, $B$ in $G$ with no common fixed points. In fact we can also have $d\neq0$ and $a\neq0$ by the same arguments.

\textbf{Case I}: $c$ is real.

Since every element of $G$ has real diagonal elements, then we get
\begin{eqnarray*}
\begin {array}{ll}
t_{4}=ta\overline{j}+b\overline{h}+c\overline{g}/t\in \mathbf{R}
\\
t_{5}=tg\overline{c}+h\overline{b}+j\overline{a}/t\in \mathbf{R}
\end{array}
\end{eqnarray*}

Solving for $b\overline{h}$ and $\overline{g}$, we find that $b\overline{h}$ and $\overline{g}$ are real.

 Next, we calculate the three cross-ratios of points $ B(0), \infty, 0,  B(\infty)$. It is easy to see that

$$\mathbb{X}_{1}=g\overline{c},\ \mathbb{X}_{2}=a\overline{j},\ \mathbb{X}_{3}=\frac{aj}{gc}.$$
Because $a,g,c,j$ are real, we know that $\mathbb{X}_{1},\mathbb{X}_{2},\mathbb{X}_{3}$ are real and $\mathbb{X}_{3}=\mathbb{X}_{2}/\mathbb{X}_{1}$.
Therefore we prove that the quadruple of point $0, \infty, B(\infty), B(0)$ lies in the same Lagrangian plane  $\mathbf{H}_{\mathbf{R}}^{2}$ from Proposition 2.2.

As
$$B(0)= \begin{bmatrix} c/j\\
f/j\\1\end{bmatrix}
, \ B(\infty)= \begin{bmatrix} a/g\\
d/g\\1\end{bmatrix}\in \mathbf{H}_{\mathbf{R}}^{2},$$
we obtain $f/j, d/g\in \mathbf{R}$. Then $f, d\in \mathbf{R}$.

Note that $$B^{2}= \begin{bmatrix} a^{2}+bd+cg& *& *\\
*&* & *\\ * & * & *\end{bmatrix},$$ we have $a^{2}+bd+cg\in\mathbf{R}$, then $b\in \mathbf{R}$. We also have $b\neq0$ by $a\overline{c}+b\overline{b}+c\overline{a}=0$.

Since $b\overline{h}$ is real,
then $h\in \mathbf{R}$. So

$$B= \begin{bmatrix} a& b& c\\
d& e& f\\g&h& j\end{bmatrix} \in \mathbf{SO}(2,1),$$
and $\langle A,B \rangle \subset \mathbf{SO}(2,1)$.

Let
$$B^{\ast}= \begin{bmatrix} a'& b'& c'\\
d'& e'& f'\\g'&h'& j'\end{bmatrix}$$
be any other element of $G$. As above, we compare the traces of $A$, $B^{\ast}$, and $AB^{\ast}$ to show that
$a',e',j'$ are real. Now
$$BB^{\ast}= \begin{bmatrix} aa'+bd'+cg'& * & *\\
*&db'+ee'+ fh'&*\\*& *&gc'+hf'+ jj'\end{bmatrix}$$
and these diagonal elements are real.
Since $a,a',e,e',j,j'$ are real, $bd'+cg', db'+ fh', gc'+hf'$ are real.

Next, we consider the element
$$B^{-1}B^{\ast}= \begin{bmatrix} ja'+fd'+cg'& * & *\\
*&hb'+ee'+ bh'&*\\*& *&gc'+df'+ aj'\end{bmatrix}.$$

It is ease to see that $fd'+cg'$, $hb'+ bh'$, $gc'+df'$ are real.
 If $fc\neq bc$, then we know that $d', g'$ are real from that $fd'+cg'$ and $bd'+cg'$ are real.
 If $fc= bc$, one can use $AB^{-1}B^{\ast}$ instead of  $B^{-1}B^{\ast}$.
 Similarly, we know that $b', c', f', h'$ are real . So $B^{\ast}$ is in $\mathbf{SO}(2,1)$. This shows that
every element of $G$ preserves $\mathbf{H}_{\mathbf{R}}^{2}$.

\textbf{Case II}: $c$ is a pure imaginary number.

We know that  $g$ is also purely imaginary and $b=d=f=h=0$. Then $$B = \begin{bmatrix} a& 0& c\\
0& e& 0\\g&0& j\end{bmatrix}$$ and $A$ leave invariant a Complex line $L$ of polar vector
\[\begin{bmatrix} 0\\1\\0
\end{bmatrix}.\]

Let
$$B_{\ast}= \begin{bmatrix} a'& b'& c'\\
d'& e'& f'\\g'&h'& j'\end{bmatrix}$$
be any other element of $G$ with real diagonal elements. Now the diagonal elements of
$$BB_{\ast}= \begin{bmatrix} aa'+cg'& * & *\\
*&db'+ee'&*\\ * & * &gc'+ jj'\end{bmatrix}$$
 are real.
So $gc'+ jj'$ is real. If $c'\neq 0$, then $c'$ is a pure imaginary number. Similarly, we have $b'=d'=f'=h'=0$.
  If $c'= 0$,  then we have $b'=f'=0$. Thus if $d'$ or $h'$ is not zero. Then $B_{\ast}$
and $A$ share exactly one common fixed point $\infty$. Therefore the subgroup $\langle B_{\ast},A\rangle$ is not discrete. So we also have $b'=d'=f'=h'=0$.
Thus we conclude that $G$ leaves invariant a Complex line $L$.
\end{pf}

 \vskip 12pt {\bf Acknowledgement}. I thank Joonhyung Kim for pointing out some mistakes in the previous  versions of the paper.
 I also thank Jieyan Wang  for  helpful discussion. This
research was supported by National Natural Science Foundational of
China (No.10671059) and  Tianyuan Foundational(No. 11126195).

\end{document}